\documentclass[english]{article}
\usepackage[T1]{fontenc}
\usepackage[latin9]{inputenc}
\usepackage{amstext}
\usepackage{amssymb}

\makeatletter

\providecommand{\tabularnewline}{\\}

\makeatother

\usepackage{babel}
\begin{document}
\title{The Biker-Hiker Problem}
\author{Peter M. Higgins}
\date{{}}
\maketitle
\begin{abstract}
There are $n$ travellers who have $k$ bicycles and they wish to
complete a journey in the shortest possible time. We investigate optimal
solutions of this problem where each traveller cycles for $\frac{k}{n}$
of the journey. Each solution is represented by an $n\times n$ binary
matrix $M$ with $k$ non-zero entries in each row and column. We
determine when such a matrix gives an optimal solution. This yields
an algorithm deciding the question of optimality of complexity $O(n^{2}\log n)$.
We introduce three symmetries of matrices that preserve optimality,
allowing identification of minimal non-optimal members of this class.
An adjustment to optimal solutions that eliminates unnecessary handovers
of cycles is established, which maintains all other features of the
solution. We identify two mutually transpose solution types, the first
uniquely minimises the number of handovers, while the second keeps
the number of separate cohorts to three while bounding their overall
separation, in the case $2k\leq n$, to under $\frac{2}{n}$ of the
journey. 
\end{abstract}

\section{The problem: not enough bicycles}

There are $n$ friends who have $k$ bicycles between them and the
group needs to reach its destination as soon as possible. How should
they go about doing this? An early allusion to this problem is in
the novel \emph{The Great House }by Cynthia Harnett {[}2{]}. Here
a pair of 17th century travelling companions with only one horse between
them adopt the `ride and tie' method for their journey from Henley-on-Thames
to London. 

\vspace{.2cm}

\textbf{Assumptions }Every person walks and cycles at the same speed
as all the others, and cycling is faster than walking. We assume that
the time required to swap from one form of locomotion to the other
is negligible. For brevity, individual travellers will sometimes be
referred to as `he' while a set of travellers will be referred to
as `they'. 

\vspace{.2cm}

\textbf{Solution }Suppose that we devise a scheme, we shall call it
an \emph{optimal }scheme, in which each traveller cycles for $\frac{k}{n}$
of the length of the journey and never stops moving forward at any
stage. Each will then have cycled and walked the same distance as
each of their companions and so all $n$ friends will arrive at their
destination simultaneously. We claim that, if it exists, such a scheme
is truly optimal in that it delivers the entire group to its destination
in the least possible time, and that any non-optimal scheme is inferior
in this respect.

It is convenient to consider the length of the journey to be $n$
units, (although we will consider divisors other than $n$). To see
that a optimal scheme is best, note that the maximum (net) forward
progress by bicycle of any scheme is $kn$. It follows that if one
member of the group of $n$ travellers cycled more than $k$ units,
then some other member must cycle less than $k$ units. This latter
traveller would then take longer than others who have cycled $k$
units (or more). Hence any approach that involved any member cycling
forward a total distance other than $k$ units would take longer to
deliver the entire group to their destination as opposed to an approach
that adopted an optimal scheme. 

That cycling is faster than walking makes the problem more interesting,
a fact that is highlighted by considering the \emph{Backpack-hiker}
problem. Here there are $k$ heavy backpacks to be transported to
the finish and any traveller carrying a backpack walks more slowly
than one that is unencumbered. The change in relative speeds makes
this problem much simpler and less interesting as for \emph{any} value
of $k$ $(1\leq k\leq n)$ the minimum time for the group to complete
the journey is the length of the journey divided by the speed of a
backpack walker.

In Section 2 we list the properties of optimal schemes more formally
through a discretised representation of the Biker-hiker problem based
on square binary matrices. In Section 3 we characterise those matrices
that correspond to optimal solutions and show that we may decide the
question of optimality for a given matrix with an algorithm that involves
$O(n^{2}\log n)$ comparisons of partial sums of the rows of the matrix.
We identify three symmetries of these optimal schemes, which leads
to the discovery of minimal schemes that assign $k$ cycled stages
to each traveller and $k$ cyclists to each stage but are nonetheless
not optimal. In Sections 4 and 5 we identify and investigate a certain
mutually transpose pair of optimal matrices for arbitrary values of
the parameters $n$ and $k$. Section 6 looks at certain facets of
these special schemes. 

\section{k-uniform solutions}

It will be convenient to allocate a measure of $n$ units for the
total length of the road the travellers will take, which we may take
to be either linear or a circuit. Along the length of the journey
we imagine there to be $n+1$ equally spaced \emph{staging posts }$P_{0},P_{1},\cdots,P_{n}$,
with $P_{0}$ and $P_{n}$ marking the beginning and end of the trip
respectively, so that the distance between successive signposts is
$1$ unit. We assign numbered symbols to each of the $n$ \emph{travellers}
as we shall call them, $t_{1},t_{1},\cdots,t_{n}$. 

\vspace{.2cm}

\textbf{Definition 2.1 }(a) The problem of delivering the $n$ travellers
equipped with $k$ bicycles $(0\leq k\leq n)$ to their common destination
in a way that minimizes the time of the last arrival will be called
the $(n,k)$-\emph{problem}.\emph{ }

(b) The leg of the journey from $P_{j-1}$ to $P_{j}$ is called \emph{stage
}$j$ and is denoted by $s_{j}$ $(1\leq j\leq n$). 

(c) An $n$-\emph{scheme $S$ }is one in which each traveller $t_{i}$
is directed to travel each stage $s_{j}$ $(1\leq j\leq n)$ either
on foot, or by bicycle. 

(d) The \emph{incidence matrix $M=M(S)$ }of an $n$-scheme $S$ is
the $n\times n$ binary matrix $M=(m_{i,j})$ $(1\leq i,j\leq n)$
where $m_{i,j}=0$ or $m_{i,j}=1$ according as traveller $t_{i}$
is directed to walk or cycle respectively stage $s_{j}$ from $P_{j-1}$
to $P_{j}$. We shall write $R_{i}$ and $C_{j}$ for the $i$th row
and $j$th column of $M$ respectively. 

(e) The \emph{scheme }$S=S(M)$ of an $n\times n$ binary matrix $M=(m_{i,j})$
is that in which traveller $t_{i}$ travels $s_{j}$ on foot or by
bicycle according as $m_{i,j}=0$ or $m_{i,j}=1$ $(1\leq i,j\leq n)$.
Note that $S(M(S))=S$ and $M(S(M))=M$. 

\vspace{.2cm}

\textbf{Definition 2.2} An $n\times n$ binary matrix $M=(m_{i,j})$
is \emph{$k$-uniform} if each row and each column contains exactly
$k$ entries equal to $1$. 

\vspace{.2cm}

\textbf{Proposition 2.3 }A scheme $S$ is optimal\emph{ }if and only
if 

(i) $M(S)$ is $k$-uniform and 

(ii) whenever a set of travellers $C$ arrives at a post $P_{j}$,
the number of cycles at $P_{j}$ is at least as great as the number
of $t_{i}\in C$ such that $m_{i,j+1}=1$. 

\vspace{.2cm}

\textbf{Proof} If $S$ is optimal then each traveller $t_{i}$ rides
exactly $k$ stages so that $R_{i}$ has exactly $k$ entries which
equal $1$. There are then $nk$ entries of $M$ equal to $1.$ If
it were not the case that each column had exactly $k$ non-zero entries,
then some column would contain more than $k$ $1's$, which is impossible
as no cycle may travel twice through the same stage. Therefore $M$
is $k$-uniform. As for Condition (ii), if it were violated then some
traveller would have to stop at some stage to wait for a bicycle to
arrive for their use. The time taken for their journey would then
exceed the optimal time unless he cycled more than $k$ stages, in
which case some other traveller would cycle fewer than $k$ stages,
and the overall time for the group to complete the journey would exceed
the optimal time. Hence if $S$ is optimal, both Conditions (i) and
(ii) must be met.

Conversely any scheme $S$ represented by a $k$-uniform matrix $M$
has exactly $k$ entries of $1$ in each row so that each traveller
is scheduled to ride $k$ stages. Condition (ii) ensures that the
progress of each traveller is never stalled by a required cycle being
unavailable upon arrival at a staging post. Therefore $S$ represents
an optimal solution. $\Box$

\vspace{.2cm}

\textbf{Definition 2.4 }We call a square $k$-uniform binary matrix
$M$ \emph{optimal }if $S(M)$ is optimal. 

\section{Optimal matrices and their symmetries}

\subsection*{Assignment mappings}

We will now introduce \emph{assignment mappings $\phi_{j}$ }for the
each stage $s_{j}$ $(1\leq j\leq n-1)$ of a scheme $S$. Suppose
$m_{i,j}=1$, meaning that $t_{i}$ cycles $s_{j}$. Then $\phi_{j}(i)=p$
conveys the information that $t_{p}$ will cycle $s_{j+1}$ on the
cycle left behind at $P_{j}$ by $t_{i}$. 

\vspace{.2cm}

\textbf{Definition 3.1 }Let $S$ denote an $n\times n$ scheme with
matrix $M=M(S)=(m_{i,j})$. A one-to-one partial mapping $\phi_{j}$
$(1\leq j\leq n-1)$ is an \emph{assignment mapping} for $S$ if

\[
\text{dom\ensuremath{\phi_{j}=\{i:m_{i,j}=1\},\,\,\text{ran\ensuremath{\phi=\{i:m_{i,j+1}=1\}.}}}}
\]

\vspace{.2cm}

The main result of this section characterises optimal schemes in terms
of the existence of a collection of assignment mappings that satisfy
two constraints. The first is the optional constraint that allows
a rider to stay on the same bike if he is required to ride two successive
stages. The second constraint ensures that $S$ is in accord with
Proposition 2.3. 

\vspace{.2cm}

\textbf{Theorem 3.2 }Let $M$ be an $n\times n$ binary matrix. Then
$S(M)$ is optimal if and only if $M$ is $k$-uniform for some $k$
$(0\leq k\leq n)$ and for each $j$, $(1\leq j\leq n-1)$ there exist
assignment mappings $\phi_{j}$ such that
\begin{equation}
\phi_{j}(i)=i\Leftrightarrow(m_{i,j}=m_{i,j+1}=1)\,\,\text{and}
\end{equation}
\begin{equation}
\sum_{l=1}^{j}m_{i',l}\leq\sum_{l=1}^{j}m_{i,l},\,\,\text{where \ensuremath{i'} denotes \ensuremath{\phi_{j}(i).}}
\end{equation}

\textbf{Proof} Suppose that $M$ is $k$-uniform and satisfies Conditions
(1) and (2). Suppose inductively that the scheme $S(M)$ has not failed
up to stage $s_{j},$ which holds when $j=1$ as $C_{1}$ has $k$
entries that equal $1$, and so travellers assigned to cycle $s_{1}$
may do so. 

Next consider stage $s_{j+1}$ from $P_{j}$ to $P_{j+1}$. For each
$i'$ such that $m_{i',j+1}=1$ there exists a unique $i$ such that
$m_{i,j}=1$ and $\phi_{j}(i)=i'$. By the inductive hypothesis, $t_{i}$
has arrived at $P_{j}$ by cycle without stalling. Condition (2) is
then exactly the requirement that ensures that this has occurred no
later than the arrival of $t_{i'}$ at $P_{j}$. Hence $S$ may continue
with $t_{i'}$ riding $s_{j+1}$ on the cycle that $t_{i}$ has ridden
on $s_{j}$. Therefore $s_{j+1}$ may be completed without stalling,
and the induction continues. The process will therefore end with $S(M)$
being fully executed without stalling, and so $S(M)$ is indeed optimal. 

Conversely, suppose that $S(M)$ is optimal. Then at stage $s_{j+1}$
$(j\geq0)$, for each $i'$ such that $m_{i',j+1}=1$, it is possible
for $t_{i'}$ to ride $s_{j+1}$ on a cycle that has been left at
$P_{j}$ by some traveller $t_{i}$. It follows that Condition (2)
is then met. This correspondence defines a partial one-to-one mapping:
\[
\phi_{j}^{-1}:\{i:m_{i,j+1}=1\}\rightarrow\{i:m_{i,j}=1\}.
\]
By uniformity, $\phi_{j}^{-1}$ is also surjective and so the partial
one-to-one mapping $\phi_{j}$ is an assignment mapping which satisfies
Condition (2). We now show that $\phi_{j}$ may be modified so that
it also satisfies Condition (1). The forward direction of the implication
in (1) follows from the definition of an assignment map, but the reverse
implication does not follow from the optimality of $S(M)$. 

Let us write $\phi$ for $\phi_{j}$ and, as before, abbreviate $\phi_{j}(i)$
to $i'$. Suppose then that $m_{i,j}=m_{i,j+1}=1$ but $i\neq i'$.
We consider the sequence $I=i,\phi(i),\phi^{2}(i),\cdots.$ If $I$
is a cycle, so that for some positive integer $p$, $\phi^{p}(i)=i$,
then it follows that $m_{t,j}=1=m_{\phi(t),j}$ for all $t=\phi^{k}(i)\,(k\geq0)$.
In this case we may modify $\phi$ (while retaining the same symbol
$\phi$ for the mapping) such that $\phi(t)=t$ for all $t=\phi^{k}(i)$,
in accord with Condition (1). Moreover, applying Condition (2) repeatedly
yields a cycle of inequalities that begins and ends with the same
sum, and so are in fact equalities, indicating that all the travellers
$t_{i},t_{\phi(i)},\cdots,t_{\phi^{p}(i)}=t_{i}$ arrive at $P_{j}$
simultaneously. The original assignment mapping $\phi$ instructed
this set of travellers to exchange bicycles in accord with the cycle
$I$. The modified mapping simply allows each traveller to remain
on the bike he is currently riding. 

Alternatively the sequence $I$ does not generate a cycle. Then by
definition of $\phi$ there exists a sequence of maximal length:
\[
i_{-r},i_{-r+1},\cdots i_{0}=i,i_{1}=\phi(i),i_{2},\cdots,i_{s-1},i_{s}
\]
\begin{equation}
\text{such that\,\,\ensuremath{\phi(i_{p})=i_{p+1},\,(-r\leq p\leq s-1)},\,\ensuremath{(r,s\geq1).}}
\end{equation}
In (3), $m_{i_{-r},j+1}=0=m_{i_{s},j}$ and $m_{t,j}=m_{t,j+1}=1$
for all $-r+1\leq t\leq s-1$. We now modify $\phi$ by putting
\begin{equation}
\phi(t)=t\,\,\forall\,-r+1\leq t\leq s-1
\end{equation}
\begin{equation}
\phi(i_{-r})=i_{s},
\end{equation}
for then Condition (2) holds trivially for $i=t$ as in (4), and (2)
also holds for (5) for $i=i_{-r},\,i'=i_{s}$ as applying Condition
(2) repeatedly for $\phi$ we have:
\[
\sum_{l=1}^{j}m_{i_{-r},l}\geq\sum_{l=1}^{j}m_{i_{-r+1},l}\geq\cdots\geq\sum_{l=1}^{j}m_{i_{s-1},l}\geq\sum_{l=1}^{j}m_{i_{s},l},
\]
which, in the notation of Theorem 3.2, provides the required inequality
concerning $i_{-r}$ and $i_{s}=\phi_{j}(i_{-r})=i_{-r}'$:
\[
\sum_{l=1}^{j}m_{i_{-r}',j}\leq\sum_{l=1}^{j}m_{i_{-r},l}.
\]

We modify $\phi$ for each such $i$, which is possible as the sequences
as in (3) that arise are pairwise disjoint as $\phi$ is one-to-one,
$i_{r}$ is not in the range of $\phi$, and $i_{s}$ is not in the
domain of $\phi$. Modifying $\phi$ as necessary for each $i$ such
that $m_{i,j}=m_{i,j+1}=1$ ensures that the partial one-to-one mapping
$\phi$ satisfies both Conditions (1) and (2), thereby completing
the proof. $\Box$ 

\vspace{.2cm}

\textbf{Definitions 3.3 }Let $M$ be a $k$-uniform matrix.

(i) For any $j$ $(1\leq j\leq n-1)$ we shall call an assignment
mapping $\phi_{j}$ \emph{optimal }if $\phi_{j}$ satisfies Conditions
(1) and (2) of Theorem 3.2. 

(ii) For any $j$ $(1\leq j\leq n-1)$ consider the partition of $X_{n}=\{1,2,\cdots,n\}$
induced by $M$ into the following four (possibly empty) disjoint
subsets:
\begin{equation}
X_{1,1}=\{i:m_{i,j}=m_{i,j+1}=1\},\,X_{1,0}=\{i:m_{i,j}=1,m_{i,j+1}=0\},
\end{equation}
\[
X_{0,1}=\{i:m_{i,j}=0,m_{i,j+1}=1\},\,X_{0,0}=\{i:m_{i,j}=m_{i,j+1}=0\}.
\]
When necessary, we write $X_{1,0}^{j}$ etc. to indicate that the
set refers to column $C_{j}$. 

An assignment mapping $\phi_{j}$ then satisfies the conditions that:
\begin{equation}
\text{dom\ensuremath{\phi_{j}=X_{1,1}\cup X_{1,0}},\,\,\,ran\ensuremath{\phi_{j}=X_{1,1}\cup X_{0,1}} }
\end{equation}
with $\phi_{j}$ acting identically on $X_{1,1}$ if $\phi_{j}$ is
optimal. 

(iii) We shall denote the $i$th row sum up to column $C_{j}$ by
$S_{i,j}$: 
\begin{equation}
S_{i,j}=\sum_{l=1}^{j}m_{i,l}\,(1\leq i,j\leq n).
\end{equation}
Suppress the second subscript $j$ by writing $S_{i}$ for $S_{i,j}$,
and form ordered sets, written in ascending order as:
\begin{equation}
\overline{X}_{1,0}=\{(i_{1},S_{i_{1}}),\cdots,(i_{p},S_{i_{p}}),\,S_{i_{1}}\leq\cdots\leq S_{i_{p},},i_{t}\in X_{1,0},\,(1\leq t\leq p)\}.
\end{equation}
\begin{equation}
\overline{X}_{0,1}=\{(j_{1},S_{j_{1}}),\cdots,(j_{p},S_{j_{p}}),\,S_{j_{1}}\leq\cdots\leq S_{j_{p}},j_{t}\in X_{0,1},\,(1\leq t\leq p)\}.
\end{equation}
To make each order unique, in the case of ties, we order by subscript
value, so if $S_{i_{1}}=S_{i_{2}}$ then $(i_{1},S_{i_{1}})<(i_{2},S_{i_{2}})$
for $\overline{X}_{1,0}$ if $i_{1}<i_{2}$, and similarly for $\overline{X}_{0,1}$.
We now meld these two lists to define a total order on $Y=\overline{X}_{1,0}\cup\overline{X}_{0,1}$.
The order $(Y,\leq)$ is equal to the order defined in (9) and (10)
when restricted to $\overline{X}_{1,0}$ and to $\overline{X}_{0,1}$
respectively. For $(i,S_{i})\in\overline{X}_{1,0}$ and $(j,S_{j})\in\overline{X}_{0,1}$
we define $(i,S_{i})<(j,S_{j})$ if $|S_{i}|\leq|S_{j}|$ and $(i,S_{i})>(j,S_{j})$
if $|S_{i}|>|S_{j}|$. In this way $\leq$ is indeed a linear order
on $Y$ as transitivity is readily checked by cases. 

\vspace{.2cm} 

\textbf{Definition 3.4 }(i) The reverse order, $(Y,\geq)$ of the
linear order $(Y,\leq)$ is the \emph{canonical order }of $Y$.

Let $A=\{a,b\}$ be a two-letter alphabet. 

(ii) The \emph{canonical word $w=a_{1}a_{2}\cdots a_{2p}\in A^{2p}$
$(a_{r}\in A,\text{1\ensuremath{\leq p\leq k}})$ }is defined by $a_{r}=a$
or $a_{r}=b$ according as the $r$th entry in the canonical order
belongs to $\overline{X}_{1,0}$ or to $\overline{X}_{0,1}$. 

(iii) For any word $w\in A^{m}$$(m\geq0)$ we write $|w|_{c}$ for
the number of instances of $c\in A$ in $w$. The \emph{length }of
$w$, denoted by $|w|$, is then $|w|=|w|_{a}+|w|_{b}$.

(iv) If $w\in A^{m}$ $(m\geq0)$ has a factorization $w=uv$, we
call $u$ a \emph{prefix }and $v$ a \emph{suffix }of $w$. 

(v) A word $w\in A^{2m}$ $(m\geq0)$ such that $|w|_{a}=|w|_{b}$
is called a \emph{Dyck word }if for every prefix $u$ of $w$, $|u|_{a}\geq|u|_{b}$. 

(vi) For $w\in A^{m}(m\geq0)$, the \emph{dual reverse word $\overline{w}$
}is formed by taking the reverse word $w^{R}$ of $w$ and interchanging
all instances of the letters $a$ and $b$.

\vspace{.2cm}

\textbf{Remark 3.5 }The set of all words of any length that satisfy
the conditions of (v) is called the \emph{Dyck language}. This is
the language of \emph{well-formed parentheses} in that replacing $a$
and $b$ by the left and right brackets `$($' and `$)$' respectively,
a Dyck word corresponds to a string of brackets that represents a
meaningful bracketing of some binary operation. For further information,
see {[}3{]}.

\vspace{.2cm}

\textbf{Proposition 3.6 }(i) The dual reverse word $\overline{w}$
of a Dyck word $w$ is also a Dyck word. 

(ii) There exists an optimal assignment mapping $\phi_{j}$ $(1\leq j\leq n-1)$
if and only if the canonical word $w=w_{j}$ is a Dyck word. 

\vspace{.2cm}

\textbf{Proof }(i) Let $\overline{w}=uv$, whence $w=\overline{\overline{w}}=\overline{v}\,\overline{u}.$
Since $w$ is a Dyck word, $|\overline{v}|_{a}\geq|\overline{v}|_{b}$,
whence $|\overline{u}|_{a}\leq|\overline{u}|_{b}$, and so $|u|_{a}\geq|u|_{b}$.
Hence $\overline{w}$ is a Dyck word. 

(ii) Suppose that $\phi=\phi_{j}$ is an optimal assignment mapping.
The action of this mapping induces a bijection from letters $a_{s}=a$
in the canonical word $w$ to letters $a_{t}=b$ in $w$, which acts,
by Condition (2) of Theorem 3.2, so that $a_{s}$ lies to the left
of $a_{t}$ in $w$. It follows that for any initial prefix $u$ of
$w=uv$, we must have $|u|_{a}\geq|u|_{b}$, for if $|u|_{a}<|u|_{b}$,
there would be some instance of $b$ in $u$ that was not in the range
of the induced mapping, contradicting that $\phi_{j}$ is one-to-one.
Hence $w$ is a Dyck word.

Conversely, given that $w$ is a Dyck word, we map $i\in X_{1,0}$
to $i'\in X_{0,1}$ whereby if $i$ corresponds to the $r$th instance
of $a$ in $w$, then $i'$ corresponds to the $r$th position of
$b$ in $w$. By the given condition, the $r$th $a$ in $w$ lies
to the left of the $r$th $b$ in $w$, whence $|S_{i}|\geq|S_{i'}|$.
The map $\phi$ thereby defined satisfies Condition (2) of Theorem
3.2. Extending $\phi$ to act identically on $X_{1,1}$ then produces
a required optimal assignment map. $\text{\ensuremath{\Box}}$ 

\vspace{.2cm}

\textbf{Theorem 3.7 }Algorithm to decide optimality of a $k$-uniform
matrix $M$. 

\vspace{.2cm}

For the columns $C_{j}$ $(1\leq j\leq n-1)$ of $M$:

1. Calculate the partial sums $S_{i,j}$ $(i\in X_{1,0}\cup X_{0,1});$

2. Rank the $2p$ $(0\leq p\leq k)$ partial sums from Step 1 in descending
order, with members of $X_{1,0}$ taking precedence over members of
$X_{0,1}$ in the case of a tie, as per Definition 3.3(iii). 

3. Form the canonical word $w=w_{j}=a_{1}\cdots a_{2p}$ where $a_{r}=a$
or $b$ according as the $r$th member of this ranking lies in $X_{1,0}$
or $X_{0,1}$. 

4. $M$ is optimal if and only if $w_{j}$ is a Dyck word for all
$1\leq j\leq n-1$. 

\vspace{.2cm}

However, it is not necessary to check the first two nor the last two
assignment mappings for optimality by virtue of part (ii) of our next
result. 

\vspace{.2cm}

\textbf{Lemma 3.8 }(i) For a given $j$ $(1\leq j\leq n-1)$, all
assignment mappings $\phi_{j}$ are optimal if and only if the canonical
word $w_{j}=a^{p}b^{p}$, $(p=|X_{1,0}|)$. 

(ii) An assignment mapping $\phi_{j}$ is optimal if $j\in\{1,2,n-2,n-1\}$
or if $k\in\{1,2,n-2,n-1\}$. 

\vspace{.2cm}

\textbf{Proof }(i) Every $\phi_{j}$ is optimal if and only if $S_{i_{1},j}\geq S_{i_{2},j}$
for all $i_{1}\in X_{1,0}$ and $i_{2}\in X_{0,1}$, which in turn
is equivalent to $w_{j}=a^{p}b^{p}$, where $p=|X_{1,0}|$. 

(ii) Let $i_{1}\in X_{1,0}$ and $i_{2}\in X_{0,1}.$ For $\phi_{1}$
and $\phi_{2}$ we have $S_{i_{1},j}\geq1$ and $S_{i_{2},j}\leq1$
$(j=1,2)$ whence it follows that $w_{j}=a^{p}b^{p}$. For $\phi_{n-2}$
or $\phi_{n-1}$ we have $S_{i_{1},j}\geq k-1$ while $S_{i_{2},j}\leq k-1$,
$(j=n-2,n-1)$ and again $w_{j}=a^{p}b^{p}.$ The claim now follows
from part (i). 

Similarly if $k\leq2$ then $S_{i_{1},j}\geq1$ and $S_{i_{2},j}\leq1$,
while if $k\geq n-2$ then $S_{i_{1},j}\geq j-1$ and $S_{i_{2},j}\leq j-1$
and again the result follows. $\Box$

\vspace{.2cm}

\textbf{Corollary 3.9 }(i) An $n\times n$ uniform matrix $M$ is
optimal if $n\le5$. 

(ii) For any non-optimal $k$-uniform matrix $M$, $3\leq k\leq n-3$. 

(iii) Optimality of a $k$-uniform matrix $M$ is preserved under
the exchange of columns $C_{1}$ and $C_{2},$ and under the exchange
of columns $C_{n-1}$ and $C_{n}$. 

\vspace{.2cm}

\textbf{Proof }(i) For $n\leq5$, for any scheme there are at most
$5-1=4$ assignment mappings which are among the four mappings listed
in Lemma 3.8(ii), and so all are optimal.

(ii) This follows from Lemma 3.8(ii).

(iii) Indeed we may replace $C_{1}$ and $C_{2}$ by any pair of binary
columns that retains $k$-uniformity of $M$, for then the transformed
matrix retains its status with respect to optimality by Lemma 3.8(ii).
These correspond to exchanging adjacent instances of $0$ and $1$
in the two columns in opposite pairs. In particular, since complete
exchange of $C_{1}$ and $C_{2}$ retains $k$-uniformity, the result
follows, as it does likewise for the exchange of the final column
pair. $\Box$

\vspace{.2cm}

\textbf{Proposition 3.10 }Let $S=S(M)$ be an $(n,k)$-uniform scheme
with a given set of assignment mappings $\phi_{j}$ $(1\leq j\leq n-1)$.
If all travellers complete $c$ cycled stages of $S$ without the
scheme failing, (that is, without any traveller being stalled) then
the scheme, with this set of assignment mappings, will not fail before
some traveller is due to ride their $(c+3)$rd cycled stage.

In particular, $S$ will not fail prior to some traveller being due
to ride their $3$rd stage, and if all travellers complete $k-2$
stages without $S$ failing, then $S$ is an optimal scheme, which
is realised by the given set of assignment mappings. 

\vspace{.2cm}

\textbf{Proof }Suppose all travellers have completed $c$ cycled stages
without failure in $S$. Suppose a walking traveller $t_{i'}$ arrives
at a staging post $P_{j}$ $(1\leq j\leq n-1)$, where $s_{j}$ represents
cycle stage number $c+1$ or $c+2$ for that traveller. Let $i=\phi_{j}^{-1}(i')$.
Then $S_{i',j}=c$ in the first case, and $S_{i',j}=c+1$ in the second.
If $t_{i'}$ stalls at $P_{j}$ then it follows that $S_{i,j}\leq c$.
However, since $m_{i,j}=1$, it follows that $t_{i}$ has not yet
completed $c$ cycled stages when the stall occurs, contrary to hypothesis.
Therefore if \emph{all} travellers complete $c$ cycle stages without
the scheme failing, then the scheme will not fail prior to some traveller
attempting to cycle a stage for the $(c+3)$rd occasion. The final
statement simply draws attention to the special cases where $c=0$,
and where $c=k-2.$ $\Box$ 

\vspace{.2cm}

\textbf{Examples 3.11 }It follows from Corollary 3.9 that the smallest
dimension $n$ that might admit a non-optimal matrix $M$ is $n=6$.
In this case, the inequality of Corollary 3.9(ii) becomes $3\leq k\leq6-3$,
so that $k=3$. Consider the simple scheme $S(M_{1})$, where $M_{1}$
is given below. This scheme is clearly optimal: travellers $t_{1},t_{2},t_{3}$
ride the first three stages and then leave their bikes to be collected
later by $t_{4},t_{5},$ and $t_{6}$ who then ride together to the
finish. The assignment mappings all act identically except for $\phi_{3}$,
which may be taken as any bijection such that $\phi_{3}(\{1,2,3\})=\{4,5,6\}$.
However, if we swap columns $C_{3}$ and $C_{4}$ in $M_{1}$, we
have the array $M_{2}$. By Lemma 3.8, the only canonical word of
$M_{2}$ that may fail to be a Dyck word is $w_{3}$. However for
$j=3$ we have $X_{1,0}=\{4,5,6\}$ and $X_{0,1}=\{1,2,3\}$. For
any $i_{1}\in X_{1,0}$ and $i_{2}\in X_{0,1}$ we have $S_{i_{1},3}=1<2=S_{i_{2},3}$
and so $w_{3}=b^{3}a^{3}$, which is not a Dyck word. Therefore $M_{2}$
is not optimal. Indeed this example shows that the class of optimal
matrices is not closed under permutation of columns, nor under the
taking of transpositions.

\vspace{.2cm}

\noindent$M_{1}=$~~%
\begin{tabular}{|c|c|c|c|c|c|c|}
\hline 
$P_{0}$ & $P_{1}$ & $P_{2}$ & $P_{3}$ & $P_{4}$ & $P_{5}$ & $P_{6}$\tabularnewline
\hline 
\hline 
$t_{1}$ & $1$ & $1$ & $1$ & $0$ & $0$ & $0$\tabularnewline
\hline 
$t_{2}$ & $1$ & $1$ & $1$ & $0$ & $0$ & $0$\tabularnewline
\hline 
$t_{3}$ & $1$ & $1$ & $1$ & $0$ & $0$ & $0$\tabularnewline
\hline 
$t_{4}$ & $0$ & 0 & $0$ & $1$ & $1$ & $1$\tabularnewline
\hline 
$t_{5}$ & $0$ & $0$ & $0$ & $1$ & $1$ & $1$\tabularnewline
\hline 
$t_{6}$ & $0$ & $0$ & 0 & $1$ & $\text{1}$ & $1$\tabularnewline
\hline 
\end{tabular}~~$M_{2}=$%
\begin{tabular}{|c|c|c|c|c|c|c|}
\hline 
$P_{0}$ & $P_{1}$ & $P_{2}$ & $P_{3}$ & $P_{4}$ & $P_{5}$ & $P_{6}$\tabularnewline
\hline 
\hline 
$t_{1}$ & $1$ & $1$ & $0$ & $1$ & $0$ & $0$\tabularnewline
\hline 
$t_{2}$ & $1$ & $1$ & $0$ & $1$ & $0$ & $0$\tabularnewline
\hline 
$t_{3}$ & $1$ & $1$ & $0$ & $1$ & $0$ & $0$\tabularnewline
\hline 
$t_{4}$ & $0$ & $0$ & $1$ & $0$ & $1$ & $1$\tabularnewline
\hline 
$t_{5}$ & $0$ & $0$ & $1$ & $0$ & $1$ & $1$\tabularnewline
\hline 
$t_{6}$ & $0$ & $0$ & $1$ & $0$ & $\text{1}$ & $1$\tabularnewline
\hline 
\end{tabular}

\vspace{.2cm}

\textbf{Theorem 3.12 }The question of whether an $n\times n$ binary
matrix $M$ is optimal may be decided by an algorithm of complexity
$O(n^{2}\log n)$. 

\vspace{.2cm}

\textbf{Proof }1. By inspecting rows and columns of $M,$ decide whether
$M$ is uniform, an operation of order $O(n^{2})$. 

If $M$ is uniform, we may decide optimality of $M$ by carrying out
the following procedure for each $j$ with $1\leq j\leq n-1$. 

2. Compute $S_{i,j+1}$ from $S_{i,j}$ for all $1\leq i\leq n-1$,
which consists of $n$ additions. This allows identification of the
sets $X_{0,0}^{j},X_{1,0}^{j},X_{0,1}^{j}$ and $X_{1,1}^{j}$. 

3. Form the two sets $\overline{X}_{1,0}^{j}$ and $\overline{X}_{0,1}^{j}$
and sort in descending order, a process which has time complexity
$O(n\ln n)$, as this is the least possible for any comparison algorithm
{[}1{]}, from which may be read the canonical word, $w_{j}$. 

4. At most $O(n)$ comparisons determine whether or not $w_{j}$ is
a Dyck word. 

For each $j$, the total complexity of steps $2$, $3,$ and $4$
is $O(n)+O(n\ln n)+O(n)=O(n\ln n)$. These steps are carried out $n-1$
times, (strictly speaking, by Lemma 3.8(ii), at most $n-5$ applications
are needed), which, including Step 1, yields an overall complexity
of $O(n^{2})+O(n^{2}\log n)=O(n^{2}\log n)$. $\Box$ 

\vspace{.2cm}

\textbf{Definition 3.13 }Let $M=(m_{i,j})$ be an $n\times n$ $k$-uniform
binary matrix. Let $S_{n}$ denote the symmetric group on $X_{n}$.
The $n\times n$ $k$-uniform matrices $M_{\pi}=(p_{i,j})$ ($\pi\in S_{n})$,$\,M_{r}=(r_{i,j})$,
and $\overline{M}=(d_{i,j})$ are defined by:

(i) $p_{i,j}=m_{\pi(i),j},$~~~(ii) $r_{i,j}=m_{i,n-j+1}$,~~~(iii)
$d_{i,j}=(m_{i,j}+1)\,\text{(mod \ensuremath{2})}$. \newline We
may denote $d_{i,j}$ by $\overline{m}_{i,j}$. 

\vspace{.2cm}

\textbf{Theorem 3.14 }Suppose that $S(M)$ is an optimal scheme. Then
so are the schemes (i) $S(M_{\pi})$, (ii) $S(M_{r})$, and (iii)
$S(\overline{M})$. 

\vspace{.2cm}

\textbf{Lemma 3.15 }Let $M$ be an $n\times n$ $k$-uniform matrix.
Then

(i) The $j$th canonical word of $M_{\pi}$ $(\pi\in S_{n})$ is $w_{j}$,
the $j$th canonical word of $M$ $(1\leq j\leq n)$. 

(i) The $j$th canonical word of $\overline{M}$ is $\overline{w}_{j}$. 

(iii) The $j$th canonical word of $M_{r}$ is $\overline{w}_{n-j}$. 

\vspace{.2cm}

\textbf{Proof }(i) The canonical words $w_{j}$ $(0\leq j\leq n-1)$
of $M$ are defined by $(Y,\leq)$ based on the partial orders as
in (9) and (10). Replacing $M$ by $M_{\pi}$, results in replacing
each of the symbols $i_{t},j_{s}$ by $\pi^{-1}(i_{t}),$$\pi^{-1}(j_{t})$
in the sets (9) and (10). Since the value of $w_{j}$ is independent
of the naming of these symbols, each canonical word $w_{j}$ is unaltered.

(ii) Write $\overline{S}_{i,j}$ for a typical partial sum of $\overline{M}$.
Since for any matrix position $(i,j)$, $\overline{S}_{i,j}=j-S_{i,j}$
the list of inequalities in (9) and (10), apart from tied sums, is
reversed when passing from $M$ to $\overline{M}$. Moreover, $i_{1}\in X_{1,0},i_{2}\in X_{0,1}$
for $M$ if and only if $i_{1}\in X_{0,1},i_{2}\in X_{1,0}$ for $\overline{M}$.
It follows from this pair of observations that the $j$th canonical
word of $\overline{M}$ is $\overline{w}_{j}$, the dual reverse canonical
word of $w_{j}$. 

(iii) Denote the partial sums of $M_{r}$ by $S_{i,j}^{r}$. Then
$S_{i,j}^{r}+S_{i,n-j}=k$ $(1\leq j\leq n,\,\text{taking\,\,}S_{i,0}=0)$.
Moreover $i_{1}\in X_{1,0},i_{2}\in X_{0,1}$ for $M$ if and only
if $i_{1}\in X_{0,1},i_{2}\in X_{1,0}$ for $M_{r}$. Now 
\[
S_{i_{1},j}^{r}\leq S_{i_{2},j}^{r}\Leftrightarrow k-S_{i_{1},n-j}\leq k-S_{i_{2},n-j}\Leftrightarrow S_{i_{2},n-j}\leq S_{i_{1},n-j}.
\]
This pair of observations imply that the $j$th canonical word of
$M^{r}$ is $\overline{w}_{n-j}$. $\Box$

\vspace{.2cm}

\textbf{Proof of Theorem 3.14. }Since $M$ is optimal, by Theorem
3.7 all canonical words $w_{j}$ of $M$ are Dyck words. By Lemma
3.15, the corresponding canonical words of $M_{\pi}$, $\overline{M}$,
and $M^{r}$ are respectively $w_{j}$, $\overline{w}_{j}$, and $\overline{w}_{n-j}$.
Since the reverse dual word of a Dyck word is a Dyck word (Proposition
3.6(i)) it follows, again by Theorem 3.7, that each of $M_{\pi}$,
$\overline{M}$, and $M^{r}$ is optimal. $\Box$ 

\vspace{.2cm}

\textbf{Definition 3.16 }Define the \emph{complementary assignment
function $\overline{\phi}_{j}$ }of an assignment function $\phi_{j}$
by putting 
\begin{equation}
\text{dom \ensuremath{\overline{\phi}_{j}=X_{0,0}\cup X_{0,1},\,\,\text{ran\ensuremath{\,\overline{\phi}_{j}=X_{00}\cup X_{1,0}}}}}
\end{equation}
with $\overline{\phi}_{j}(i)=i$ if $i\in X_{0,0}$ and $\overline{\phi}_{j}(i)=\phi_{j}^{-1}(i)$
if $i\in X_{0,1}$. 

\vspace{.2cm}

\textbf{Remark 3.17 }We may prove Theorem 3.14 directly by identifying
optimal assignment mappings $\psi_{j}$ for the matrix of the transformed
scheme in terms of given optimal assignment mappings $\phi_{j}$ of
$M(S)$. In case (iii) for instance, put $\psi_{j}=\overline{\phi}_{j}$
$(1\leq j\leq n-1)$, as per Definition 3.16. For $\overline{M}$
we have 
\[
\text{dom\ensuremath{\psi_{j}=\{i:m_{i,j}=0\}=\{i:d_{i,j}=1\},}}
\]
\[
\text{\ensuremath{\text{ran\ensuremath{\psi_{j}=\{i:m_{i,j+1}=0\}=\{i:d_{i,j+1}=1\}}},}}
\]
whence it follows that the $\psi_{j}$ qualify as assignment mappings
for $S(\overline{M})$. Moreover, by definition, $\psi_{j}(i)=i$
if and only if $d_{i,j}=d_{i,j+1}=1$, and so Condition (1) is satisfied.
For $i\in X_{0,0}$ we have $\psi_{j}(i)=i$ and so in this case the
inequality of Condition (2) becomes an equality, and is thus satisfied.
Otherwise $i\in X_{0,1}$. Then we have
\[
\sum_{l=1}^{j}d_{\psi_{j}(i),l}=j-\sum_{l=1}^{j}m_{\psi_{j}(i),l}=j-\sum_{l=1}^{j}m_{\overline{\phi}_{j}(i),l}=j-\sum_{l=1}^{j}m_{\phi_{j}^{-1}(i),l}
\]
\[
\leq j-\sum_{l=1}^{j}m_{i,l}=\sum_{l=1}^{j}d_{i,l},
\]
where the inequality comes from Condition (2) applied to the $\phi_{j}$,
thereby verifying Condition (2) for the $\psi_{j}$. For parts (i)
and (ii) the corresponding assignment mappings are given respectively
by $\psi_{j}=\pi^{-1}\phi_{j}\pi$ , and $\psi_{j}=\phi_{n-j}^{-1}$,
$(1\leq j\leq n-1)$. 

\subsection*{Removing unnecessary handovers from an optimal scheme }

Optimal schemes may have unnecessary cycle handovers, which can be
removed, resulting in a scheme that is still optimal and displays
the same character as the original. Suppose that $S=S(M)$ is an optimal
$(n,k)$-scheme and for some $j$ we have $i_{1}\in X_{1,0},i_{2}\in X_{0,1}$
and $S_{i_{1},j}=S_{i_{2},j}$. Then $t_{i_{1}}$ and $t_{i_{2}}$
arrive at $P_{j}$ simultaneously, the former by bike and the latter
on foot, whereupon $t_{i_{2}}$ takes one of the bikes parked at $P_{j}$
and goes on to cycle $s_{j+1}$. However, one cycle handover could
be avoided if the pair of travellers swapped labels at this point,
with $t_{i_{1}}$ taking on the mantle of $t_{i_{2}}$ and vice-versa.
In other words $t_{i_{1}}$ would complete the journey as instructed
by the final part of $R_{i_{2}}$ from $m_{i_{2},j+1}$ onwards and
similarly $t_{i_{2}}$ would follow $R_{i_{1}}$ from $m_{i_{1},j+1}$
onwards, allowing $t_{i_{1}}$ to remain on his bike for $s_{j+1}$. 

This does not alter any column sums, and nor does it alter rows sums
as the initial portions are equal: $S_{i_{1},j}=S_{i_{2},j}$, and
hence so are the latter portions, as together they each sum to $k$.
Applying this procedure repeatedly will lead to a more efficient scheme
that will \emph{appear} to be identical, meaning that if both schemes
were to run simultaneously, at any given moment the set of positions
of walking travellers and the set of positions of cycling travellers
for the two schemes are identical. We shall call such a scheme \emph{reduced},
with it being free of \emph{excess handovers}. In summary we have
the following theorem. 

\vspace{.2cm}

\textbf{Theorem 3.18 }Given any optimal scheme $S=S(M)$ for the $(n,k)$-problem
we may construct an optimal scheme $S(M')$ that is free of unnecessary
handovers by repetition of the rule that if for some $j$ we have
$i_{1}\in X_{1,0},i_{2}\in X_{0,1}$ and $S_{i_{1},j}=S_{i_{2},j}$
we replace $R_{i_{1}}$ and $R_{i_{2}}$ in $M$ by
\[
R_{i_{1}}^{'}=(m_{i_{1},1},\cdots,m_{i_{1},j},m_{i_{2},j+1},\cdots,m_{i_{2},n}),
\]
\begin{equation}
R_{i_{2}}^{'}=(m_{i_{2},1},\cdots,m_{i_{2},j},m_{i_{1},j+1},\cdots,m_{i_{1},n}).
\end{equation}

\textbf{Remark 3.19} Removal of unnecessary handovers yields a stronger
form of Condition (2) of Theorem 3.2 in which all the associated inequalities
for which $i'\neq i$ are strict, for \emph{all }collections of optimal
assignment mappings. However, this process does alter the scheme,
whereas imposing Condition (1) merely chooses a special type of set
of assignment maps for a given scheme.

Conversely, if $S(M)$ is optimal and \emph{every} set of optimal
assignment mappings yields strict inequalities in Condition (2), it
follows that $S(M)$ has no unnecessary handover. However an optimal
scheme may have \emph{some }collection of assignment mappings for
which the non-trivial inequalities in Condition (2) are all strict,
yet the scheme still not be reduced. Such a collection of assignment
mappings has the added feature that each traveller will find a parked
cycle waiting for him whenever he is due to pick one up. 

Simple camparison arguments like those in the proof of Theorem 3.14
give the following result. 

\vspace{.2cm}

\textbf{Proposition 3.20 }For any optimal matrix $M$, the number
$h=h(M)$ of excess handovers is the same for the optimal schemes
$M_{\pi}$, $M_{r}$ and $\overline{M}$.

\section{Solution to the Biker-hiker Problem}

We now provide a particular solution type to the general Biker-hiker
problem. Because of the cyclic nature of our solutions, it will be
convenient in this section to label the travellers as $t_{0},t_{1},\cdots,t_{n-1}$
and the entries of an $n\times n$ matrix $M$ as $m_{i,j}$ $(0\leq i,j\leq n-1)$,
and stages are labelled $s_{0},s_{1},\cdots,s_{n-1}$ also. 

\vspace{.2cm}

\textbf{Definition 4.1 The Cyclic Scheme} We define the cyclic $(n,k)$-scheme
$S=S_{n,k}$ with matrix $M(S)=M_{n,k}$ by assigning the cycling
quota for $t_{i}$ to consist of the $k$ cyclically successive stages,
which run from $P_{ik}$ to $P_{(i+1)k}$, where arithmetic is conducted
modulo $n$. The matrix $M_{n,k}$ of the $n\times n$ cyclic scheme
will be called the \emph{cyclic} $(n,k)$\emph{-matrix}. 

\vspace{.2cm}

Since we are working modulo $n$, we identify $P_{0}$ and $P_{n}$,
thereby making the journey a circuit. However, the following analysis
holds whether the journey is linear or circular in nature. 

\vspace{.2cm}

\textbf{Theorem 4.2 }The $(n,k)$-cyclic scheme $S_{n,k}$ is optimal. 

\vspace{.2cm}

\textbf{Proof} By construction, $M=M_{n,k}$ is row $k$-uniform.
The entry $m_{i,j}=1$ if and only if $j$ belongs the cyclic sequence
$ik,ik+1,\cdots,ik+k-1$ which is equivalent to the statement that
$ki$ (mod $n$) lies in the cyclic interval $I_{j}=(j-k+1,j-k+2,\cdots,j)$.
Therefore the number of $1's$ in $C_{j}$ is the number of solutions
to the congruences $kx\equiv a$ (mod $n$), $a\in I_{j}$. Such a
congruence has no solution if $d=$ gcd$(n,k)$ is not a divisor of
$a$, otherwise there are $d$ solutions. Since $d|n$, it follows
that the number of $a$ such that $d|a$ is the number of multiples
of $d$ in $I_{j}$ when $I_{j}$ is regarded as an interval of $k$
consecutive integers, which is $\frac{k}{d}$, and so that there are
exactly $d\cdot\frac{k}{d}=k$ non-zero entries in each column of
$M$. (Indeed every column of $M$ represents the same cyclic sequence:
see Prop. 4.9.)

To prove optimality of the matrix $M$ of a cyclic scheme we appeal
to Proposition 3.10, which says that a uniform scheme will not stall
prior to some traveller attempting to mount a bicycle for the third
time. Since no-one mounts a bike more than twice in a cyclic scheme,
it follows that there is no stalling and the scheme is optimal. $\Box$

\vspace{.2cm}

For $M=(m_{i,j})$, a square matrix, $M_{r}$, the matrix that results
from reversing the rows of $M$ is described by permuting the columns
of $M$ by $C_{j}\leftrightarrow C_{n-j-1}$. Similarly we now define
$M_{c}$ by reversing the columns of $M$, which is effected by the
row permutation whereby $R_{i}\leftrightarrow R_{n-i-1}$. Of course
both these permutations are respectively involutions on the set of
columns and the set of rows of $M$. Writing $M_{rc}$ for $(M_{r})_{c}$,
and similarly defining $M_{cr},M_{r^{2}}$ and so on, we see that
$M_{rc}=M_{cr}=(m_{n-1-i,n-1-j})$. 

\vspace{.2cm}

\textbf{Proposition 4.3 }(i) For the cyclic $(n,k)$-matrix $M$,

(i) $M_{c}=M_{r}$. (ii) $M_{cr}=M_{rc}=M$. (iii) $(M^{T})_{r}=(M_{r})^{T}$,
$(M^{T})_{c}=(M_{c})^{T}$ (iv) $(M^{T})_{rc}=M^{T}$. 

\vspace{.2cm}

\textbf{Proof }We prove (i), from which (ii), (iii), and (iv) readily
follow. For $M=(m_{i,j})$ we have $M_{c}=(c_{i,j})$ where $c_{i,j}=m_{n-1-i,j}$
and $M_{r}=(r_{i,j})$, where $r_{i,j}=m_{i,n-1-j}$ . Then we have
\[
c_{i,j}=1\Leftrightarrow m_{n-1-i,j}=1\Leftrightarrow j\equiv(n-1-i)k+a\,\,\text{(mod \ensuremath{n)\,\,\text{for some \ensuremath{0\leq a\leq k-1}}}}
\]
\begin{equation}
\Leftrightarrow j+ik+k\equiv a\,\,(\text{mod \ensuremath{n)}}
\end{equation}
\[
r_{i,j}=1\Leftrightarrow m_{i,n-j-1}=1\Leftrightarrow n-j-1\equiv ik+b\,\text{(mod\,\ensuremath{n)\,\,\text{for some \ensuremath{0\leq b\leq k-1}}}}
\]
\[
\Leftrightarrow j+ik+1\equiv-b\,\,(\text{mod\,\ensuremath{n)}}
\]
\begin{equation}
\Leftrightarrow j+ik+k\equiv c\,\,(\text{mod\,\ensuremath{n),}}
\end{equation}
where $c=k-1-b$. Now 
\[
0\leq b\leq k-1\Leftrightarrow-k+1\leq-b\leq0\Leftrightarrow0\leq c\leq k-1.
\]
We now note that the conditions of (13) and (14) are the same. It
follows that $c_{i,j}=r_{i,j}$, allowing us to conclude that $M_{c}=M_{r}$.
$\Box$

\vspace{.2cm}

\textbf{Theorem 4.4 }For the $(n,k)$-problem, $(1\leq k\leq n-1)$
on an $n$-circuit, the cyclic scheme matrix represents the unique
solution, up to permutation of rows, in which each traveller mounts
and dismounts a cycle only once. 

\vspace{.2cm}

\textbf{Proof} By construction $S(M_{n,k})$ instructs each traveller
to mount and dismount a cycle exactly once on the circuit. On the
other hand, a uniform scheme that has this property is the cyclic
solution. To see this, take any traveller, label the traveller $t_{0}$
and label the post where $t_{0}$ mounts a cycle as $P_{0}$. Since
$t_{0}$ has a single bike ride, he must pass posts that we may label,
$P_{1},P_{2},\cdots$ until he alights at a post that we may label
$P_{k}$, thereby completing his full quota. That bicycle is then
picked up by another traveller, who we may label $t_{1}$, who rides
between posts that we may label $P_{k}$ to $P_{2k}$ (subscripts
modulo $n$). We continue this process with the traveller labelled
$t_{i}$ riding the $k$ stages from $P_{ik}$ to $P_{(i+1)k}$. But
this is just the description of the cyclic solution of the $(n,k$)-problem.
$\Box$

\vspace{.2cm}

\textbf{Remark 4.5 }The feature of one ride per ciruit is preserved
by any of the symmetries of Theorem 3.14. In the case of row reversal,
the non-zero stages for $t_{i}$ remain those between $P_{ik}$ and
$P_{(i+1)k}$ but are now ridden in reverse. Indeed since, by Proposition
4.3, $M_{r}=M_{c}$, we see that for the cyclic solution matrix $M$,
$M_{r}$ is a special case of permutation of the rows of $M$, and
so $M_{r}$ also represents an $(n,k)$-cyclic scheme. When we pass
to the binary dual we find that $\overline{M}_{n,k}=(M_{n,n-k})_{r}$
and so by the previous observation it follows that $\overline{M}_{n,k}$
indeed represents a cylic solution to the $(n,n-k)$ problem. In detail,
write $(M_{n,n-k})_{r}=(a_{i,j})$ and $M_{n,n-k}=(m_{i,j})$ whence
$a_{i,j}=1$ becomes
\[
m_{i,n-1-j}=1\Leftrightarrow n-1-j\equiv i(n-k)+a\,\,\text{(mod \ensuremath{n}) for some \ensuremath{0\leq a\leq n-k-1}}
\]
\begin{equation}
\Leftrightarrow j+1+a\equiv ik\,\,\text{(mod\,\ensuremath{n)\,\,\text{\ensuremath{0\leq a\leq n-k-1.}}}}
\end{equation}
For the left hand side we write $\overline{M}_{n,k}=(b_{i,j})$ and
$M=(m_{i,j})$. Then $b_{i,j}=1$ may be written as
\[
m_{i,j}=0\Leftrightarrow j\equiv(i+1)k+b\,\,\text{(mod\,\ensuremath{n)\text{ for some \ensuremath{0\leq b\leq n-k-1}}}}
\]
\[
\Leftrightarrow j-(b+k)\equiv j+(n-b-k)\equiv ik\,\,\text{(mod\,\ensuremath{n)}}.
\]
Now $0\leq n-1-b-k\leq n-1-k$. Put $c=n-1-b-k$. Then
\begin{equation}
j+1+c\equiv ik\,\,\text{(mod\,\ensuremath{n})\,\,\ensuremath{0\leq c\leq n-k-1.}}
\end{equation}
The agreement of (15) and (16) allow us to conclude that $\overline{M}_{n,k}=(M_{n,n-k})_{r}$
and so $\overline{M}_{n,k}$ represents a cyclic solution to the $(n,n-k)$-problem. 

\vspace{.2cm}

\textbf{Proposition 4.6} Consider the $(n,k)$-problem and let $d=$
gcd$(n,k)$. Let $R_{i}$ and $C_{j}$ denote the $i$th row and $j$th
column respectively of $M$, the matrix of the cyclic solution to
the $(n,k$)-problem as defined in 4.1. Then 

(i) $R_{i}=R_{j}$ if and only if $i\equiv j$ (mod $\frac{n}{d}$);

(ii) $C_{i}=C_{j}\,$ if and only if $dq\leq i,j\leq d(q+1)-1$ for
some $q\in\{0,1,\cdots,\frac{n}{d}-1\}$. 

\vspace{.2cm}

\textbf{Proof }(i) is trivially true if $k=0$ or $k=n$. Otherwise
the cyclic intervals of entries that equal $1$ in $R_{i}$ and $R_{j}$
respectively are defined by the corresponding cyclic lists of staging
posts: $P_{ik},P_{ik+1},\cdots,P_{(i+1)k}$ and $P_{jk},P_{jk+1},\cdots,P_{(j+1)k}$.
These lists are identical if and only if $ik\equiv jk$ (mod $n)\Leftrightarrow i\equiv j$
(mod $\frac{n}{d})$. 

(ii) We observe that the non-zero entries of each row $R_{i}$ consist
of two intervals: an initial interval $I$ of $R_{i}$ of length $r$
say, and a terminal interval $T$ of $R_{i}$ of length $k-r$ $(0\leq r\leq k)$.
We may write $k=du$ and $n=dv$. Then for some $x\geq0$ we have
\[
ik\,\text{(mod \ensuremath{n)=dui-dvx=d(ui-xv).}}
\]
If non-empty, the terminal interval $T$ begins at $P_{ik}$ and ends
at $P_{n}$ and so has length $|T|$ given by
\[
|T|=n-ik\text{(mod \ensuremath{n)} = \ensuremath{d(v-ui+xv)}. }
\]
It follows that $d||T|$. The length $|I|$ of the initial interval
is $|I|=k-|T|=du-|T|$, whence $d||I|$ also. In the case where both
$I$ and $T$ are non-empty the (successive) zeros in $R_{i}$ number
$n-|I|-|T|$, which likewise is a multiple of $d$. Otherwise there
is an initial interval of zeros of length $ik$, which is a multiple
of $d$, from which it follows that the terminal interval of zeros
has length that is too a multiple of $d$. Therefore within any row,
counting left to right by columns, the entries from one multiple of
$d$ up to but not including the next, are equal, because each maximal
list of identical entries begins at a multiple of $d$. Hence 

\begin{equation}
dq\leq i,j\leq d(q+1)-1\,\,(0\leq q\leq\frac{n}{d}-1)\Rightarrow C_{i}=C_{j}.
\end{equation}

In order to prove the reverse implication, we introduce the following
construction. By (17), the columns of $M$ consist of $\frac{n}{d}$
blocks $A_{1},A_{2},\cdots,A_{\frac{n}{d}}$ of contiguous columns,
with each $A_{i}$ consisting of $d$ identical columns. On the other
hand $M$ is partitioned into $\frac{n}{d}$ sets of $d$ (non-contiguous)
identical rows $B_{1},B_{2},\cdots,B_{\frac{n}{d}}$. We may permute
the rows of $M$, giving a new optimal matrix $M'$ in which the rows
of $M'$ are partitioned into $\frac{n}{d}$ blocks $B_{1}',B_{2}',\cdots,B_{\frac{n}{d}}'$
each consisting of $d$ identical rows. The new column blocks, $A_{1}',A_{2}',\cdots,A_{\frac{n}{d}}'$
that result from this row permutation each consist of $d$ columns,
and the columns within each block remain identical. The pairwise intersections
$A_{i}'\cap B_{j}'$ partition $M'$ into $\frac{n^{2}}{d^{2}}$ square
blocks, which are themselves $d\times d$ matrices. Each such block
has identical columns and identical rows, whence it follows that all
entries of any particular $A_{i}'\cap B_{j}'$ are identical. We can
then form a \emph{quotient matrix},\emph{ }$M_{d}'$ by identifying
each of the $A_{i}'\cap B_{j}'$ with the common value ($0$ or $1$)
of all entries in that sub-matrix. Therefore $M_{d}'$ is the cyclic
scheme for the $(\frac{n}{d}\times\frac{k}{d})$-problem in which
the travellers and the bicycles are grouped into sets of order $d$,
which move together as a block throughout the scheme. 

If now the reverse implication in (17) were false, it would imply
that there were two identical columns in the quotient matrix $M_{d}'$.
It is possible to prove directly by analysing the cardinality of the
intersection of sets of cylic intervals that in the case where $n$
and $d$ are coprime, no two columns are identical, which, since $(\frac{n}{d},\frac{k}{d})$
is a pair of coprime integers, applies to $M_{d'}$. However the desired
result follows at once from the next proposition which shows that
in the case of coprimality the deteminant of $M$ corresponds to the
number of bicycles. 

\vspace{.2cm}

\textbf{Proposition 4.7 }If $n$ and $k$ are coprime then $|\text{det\ensuremath{(M_{n,k})|=k}}.$
Otherwise $M_{n,k}$ is singular. 

\vspace{.2cm}

\textbf{Proof }Let $d=$ gcd$(k,n)$. If $d\geq2$ then by Proposition
4.6(i), $M_{n,k}$ has a pair of identical rows and so det$(M)=0$.
For $d=1$ however the rows are cyclically identical and no two are
equal. It follows that the set of rows consists of all $n$ different
possibilities that arise from the cyclic sequence $(1,1,\cdots,1,0,\cdots,0)$,
where the initial sequence of $1$'s has length $k$. By permuting
the rows of $M_{n,k}$ we may obtain the circulant matrix $C_{n,k}$,
where $R_{i}(C_{n,k})$ has for its non-zero entries $m_{,i,i},m_{i,i+1},\cdots,m_{i,i+k-1}$,
(addition modulo $n$). Hence det$(M_{n.k})=\pm$ det$(C_{n,k})$.
We may therefore complete the proof by showing that det$(C_{n,k})=k$.

By a standard result on circulant matrices (see for example {[}4{]}),
with $\omega$ denoting any primitive $n$th root of unity:
\begin{equation}
\text{det\ensuremath{(C_{n,k})=\Pi_{i=0}^{n-1}(1+\omega^{i}+\omega^{2i}+\cdots+\omega^{(k-1)i}).}}
\end{equation}
For $i=0$, the bracketed term is equal to $k$. It remains to show
that the product of the other terms in (18) is equal to $1$. By summing
each of the geometric series we see that this claim is equivalent
to the equation:
\begin{equation}
\text{\ensuremath{\Pi_{i=1}^{n-1}(\omega^{ki}-1)=\Pi_{i=1}^{k}(\omega^{i}-1).}}
\end{equation}
However, since $k$ and $n$ are coprime, $\omega^{k}$ is also a
primitive $n$th root of unity, and so it follows that the products
in (19) are identical up to the order of their factors, thereby completing
the proof. In particular, no two columns of $M_{n,k}$ are identical,
thereby also completing the proof of Proposition 4.6. $\Box$

\vspace{.2cm}

\textbf{Remark 4.8 }Note from the previous proof that for gcd$(n,k)=1$,
$S(C_{n,k})$ is also the cyclic $(n,k)$-scheme. Moreover, the non-zero
entries of $R_{i}(C_{n,k}^{T})$ are $m_{i,i},m_{i,i-1},\cdots,m_{i,i-k+1}$.
Hence the non-zero entries of $R_{i+k-1}(C_{n,k}^{T})$ are $m_{i+k-1,i+k-1},m_{i+k-1,i+k-2},\cdots,m_{i+k-1,i}$,
which match those of $R_{i}(C_{n,k})$, and so $S(C_{n,k}^{T})$ is
also the cyclic $(n,k)$-scheme, with $C_{n,k}^{T}$ obtained by rotating
the columns of $C_{n,k}$ forward by $k-1$ places. This contrasts
with $S(M_{n,k}^{T})$, the subject of Section 5, which although optimal
is of a different character to $S(M_{n,k})$. 

\vspace{.2cm}

The rows of $M_{n,k}$ represent the same cyclic sequence. The same
is true of the columns. 

\vspace{.2cm}

\textbf{Proposition 4.9 }Let $M=M_{n,k}$ be the cyclic $(n,k)$-matrix.
Then every pair or columns of $M$ represent the same cyclic sequence. 

\vspace{.2cm}

\textbf{Proof} Let gcd$(n,k)=d$. For $M_{n,k}=(m_{i,j})$ we have,
with addition modulo $n$, that $m_{i,j}=m_{i+1,j+k}$. Since gcd$(n,k)=d$,
there exists a value $r$ such that $kr\equiv d$ (mod $n)$; $r$-fold
application of the previous equation then gives $m_{i,j}=m_{i+r,j+kr}=m_{i+r,j+d}$.
It follows in particular that $C_{j}$ and $C_{j+d}$ define the same
cyclic sequence, with one being transformed into the other through
a rotation of $r$ positions. By Proposition 4.7(ii), the columns
$C_{0},C_{1},\cdots,C_{d-1}$ are identical, and so it now follows
that every pair of columns of $M_{n,k}$ define the same cyclic sequence.
$\Box$

\section{The transpose solution}

We have noted that optimality of a uniform matrix is generally not
preserved under transposition. However, the cyclic scheme is an exception
to this.

\vspace{.2cm}

\textbf{Theorem 5.1 }The transpose matrix $M=M_{n,k}^{T}$ of a cyclic
$(n,k)$-matrix $M_{n,k}$ is also optimal. 

\vspace{.2cm}

We shall call $S(M^{T})$ a \emph{transpose cyclic scheme }and similarly
$M^{T}$ is a \emph{transpose cyclic matrix}. With subscripts calculated
modulo $n$, the non-zero entries of column $C_{j}$ of $M^{T}$ are
$m_{jk,j},m_{jk+1,j},\cdots,m_{jk+k-1,j}$ $(0\leq j\leq n-1)$. The
transpose cyclic matrix $M^{T}$ is $k$-uniform, and so if $S(M^{T})$
does not stall, we have optimality. By passing to the binary dual
if necessary, we may suppose that $k\leq\frac{n}{2}$, for first note
that for any binary matrix $M=(m_{i,j}),$ we have $\overline{M}^{T}=\overline{M^{T}}$
as the $(i,j)$th entry in each of these matrices is $\overline{m}_{j,i}$.
Now let us assume that for any cyclic $(n,k)$-matrix $M$ with $n\geq2k$,
the transpose matrix $M^{T}$ is optimal. Suppose that $M$ is a cyclic
$(n,k)$-matrix with $n<2k$ and consider $M^{T}$. Then $\overline{M^{T}}=\overline{M}^{T}$,
with $\overline{M}$ a cyclic $(n,n-k)$-matrix. Since $n<2k$, it
follows that $n>2(n-k)$, and so by our assumption we have that $\overline{M}^{T}$
is optimal. But $\overline{M}^{T}=\overline{M^{T}}$, whence $\overline{\overline{M^{T}}}=M^{T}$
is also optimal. 

We are therefore permitted to adopt the assumption that $2k\leq n$
in our proof that transpose cyclic matrices are optimal. For the remainder
of the section we shall denote our transpose cyclic matrix by $M$
(as opposed to $M^{T}$). For any $t$, at least one of the entries
$m_{t,j}$ and $m_{t,j+1}$ of $M$ is $0$, as we now show. 

For any $j\geq0$, there is a unique $i\,(=jk$ mod $n)$, such that
the non-zero entries of columns $C_{j}$ and $C_{j+1}$ in $M$ have
the form:

\[
(m_{i,j}=m_{i+1,j}=\cdots=m_{i+k-1,j}=1)
\]
\begin{equation}
\Leftrightarrow(m_{i+k,j+1}=m_{i+k+1,j+1}=\cdots=m_{i+2k-1,j+1}=1).
\end{equation}
Since the total number of entries listed in (20) is $2k\leq n$, it
follows that there is no $t$ such that $m_{t,j}=m_{t,j+1}=1$, as
claimed. 

The non-zero entries of $C_{j}$ form a cyclic block of length $k$.
This will manifest itself either as a single linear block in $C_{j}$,
or as a pair of \emph{initial }and a \emph{terminal }blocks. In the
single block case, the initial and terminal blocks of non-zero entries
are one and the same. 

\vspace{.2cm}

\textbf{Lemma 5.2 }Let $(i,j)$ be the final entry of the initial
block of non-zero entries of $C_{j}$. We shall write $i=i(j)$. Then
\begin{equation}
S_{0,j}=S_{1,j}=\cdots=S_{i,j}=S_{i+1,j}+1;\,S_{i+1,j}=S_{i+2,j}=\cdots=S_{n-1,j}.
\end{equation}

\textbf{Proof }We proceed by induction on $j$. For $j=0$ we have
$m_{0,0}=m_{1,0}=\cdots=m_{k-1,0}=1$, $m_{k,0}=\cdots=m_{n-1,0}=0$,
in accord with (21), where $i(0)=k-1$. Suppose now that (21) holds
for some value of $j$ and consider $C_{j+1}$. Suppose first that
the non-zero entries of $C_{j+1}$ form a single linear block: $m_{t,j+1}=m_{t+1,j+1}=\cdots=m_{t+k-1,j+1}=1$.
If $t=0$ then $i(j)=n-1$ in (21) and all the row sums for $C_{j}$
in (21) are equal. It then follows that (21) holds for $C_{j+1}$
as in the $j=0$ case. Otherwise $t\geq1$ and so $i(j)=t-1$. By
induction:
\[
S_{0,j}=S_{1,j}=\cdots=S_{t-1,j}=S_{t,j}+1,\,\,S_{t,j}=S_{t+1,j}=\cdots=S_{n-1,j}.
\]
Since $S_{p,j}=S_{p,j+1}$ for all $0\leq p\leq t-1$ it follows that
$S_{0,j+1}=S_{p,j+1}$ for all $0\leq p\leq t-1$. On the other hand
for $t\leq p\leq t+k-1$ we have $S_{p,j+1}=1+S_{p,j}=1+(S_{0,j}-1)=S_{0,j}=S_{0,j+1}$.
Therefore $S_{0,j+1}=S_{p,j+1}$ for all $0\leq p\leq t+k-1$. Finally,
for the case where $t+k\leq p$ we have $S_{p,j+1}=S_{p,j}=S_{0,j+1}-1$
and so (21) is holds for the $S_{p,j+1}$ $(0\leq p\leq n-1)$.

The alternative case is where the non-zero entries of $C_{j+1}$ break
into distinct initial and terminal blocks. The two blocks then have
the respective forms: 

\[
m_{0,j+1}=m_{1,j+1}=\cdots=m_{i,j+1}=1
\]
\begin{equation}
\&\,\,m_{n-k+i+1,j+1}=m_{n-k+i+2,j+1}=\cdots=m_{n-1,j+1}=1\,\,(0\leq i\leq k-2).
\end{equation}
(Note that the total number of entries in (22) is indeed $(i+1)+(n-1-(n-k+i))=k$.)
The single linear cyclic block of non-zero entries of $C_{j}$ ends
at $m_{n-k+i,j}=1$ and begins at $m_{(n-k+i-(k-1)),j}=m_{n-2k+i+1,j}=1$.
By applying the inductive hypothesis to the row sums of $C_{j}$ we
infer that:
\begin{equation}
S_{0,j+1}=S_{1,j+1}=\cdots=S_{i,j+1}=S_{0,j}+1=\cdots=S_{i,j}+1.
\end{equation}
Since $S_{0,j}=S_{1,j}=\cdots=S_{n-k+i,j}$, it follows that
\begin{equation}
S_{i+1,j+1}=\cdots=S_{n-k+i,j+1}=S_{0,j}.
\end{equation}
Finally we have 
\[
S_{n-k+i,j}=1+S_{n-k+i+1,j},\,\,S_{n-k+i+1,j}=\cdots=S_{n-1,j}
\]
\begin{equation}
\Rightarrow S_{n-k+i+1,j+1}=\cdots=S_{n-1,j+1}=1+S_{n-k+i+1,j}=S_{0,j}.
\end{equation}
Statements (23), (24), and (25) together give (22) as applied to $C_{j+1}$. 

Now suppose that $n\leq2k$. Recall from Remark 4.5 that the binary
dual $\overline{M}$ of $M$ is the cyclic transpose matrix of the
$(n,n-k)$ problem with rows reversed. 

Since $n\geq2(n-k)$ it follows that (22) holds for the corresponding
row sums of the columns of $\overline{M}$, (denoted $\overline{S}_{i,j}$).
Note that $S_{i,j}+\overline{S}_{i,j}=j+1$. Hence for some $i$ $(0\leq i\leq n-1)$:
\[
\overline{S}_{n-1,j}=\overline{S}_{n-2,j}=\cdots=\overline{S}_{n-i,j}=\overline{S}_{n-i-1,j}+1;\,\overline{S}_{n-i-1,j}=\overline{S}_{n-i-2,j}=\cdots=\overline{S}_{0,j},
\]
\[
\Leftrightarrow S_{n-1,j}=S_{n-2,j}=\cdots=S_{n-i,j}=S_{n-i-1,j}-1;
\]
\[
S_{n-i-1,j}=S_{n-i-2,j}\cdots=S_{0,j},
\]
\[
\Leftrightarrow S_{0,j}=\cdots=S_{n-i-1,j}=S_{n-i,j}+1;\,S_{n-i,j}=\cdots=S_{n-1,j},
\]
which is in accord with (22) with $i(j)=n-1-i$. This completes the
proof. $\Box$ 

\vspace{.2cm}

\textbf{Proof of Theorem 5.1} As already observed, we may assume that
$n\geq2k$, in which case it is clear from (21) that for any column
$C_{j},$ 
\[
(i_{1}\in X_{1,0}^{j},i_{2}\in X_{0,1}^{j})\Rightarrow S_{i_{1},j}\geq S_{i_{2},j},
\]
from which it follows that the canonical word $w_{j}$ is the Dyck
word $w_{j}=a^{k}b^{k}$ $(2k\leq n\Rightarrow|X_{1,0}|=k)$. Hence
every assignment mapping $\phi_{j}$ is optimal, and therefore $M$
is optimal. $\Box$

\vspace{.2cm}

\textbf{Proposition 5.3 }If $n\geq2k$, then at any time point during
the execution of $S(M)$, there are at most 3 distinct positions for
the travellers. Moreover the distance separating one cohort from the
next is less than $1$ unit. 

\vspace{.2cm}

\textbf{Proof }We begin with three useful observations.

$\bullet$ If $i\leq i'$ then $t_{i}$ never trails $t_{i'}$. This
follows easily from the fact that for any fixed $j$, the $S_{i,j}$
are monotonically decreasing in $i$ (Lemma 5.2).

$\bullet$ Consider a typical column $C_{j}$ of $M$. As explained
prior to Lemma 5.2, $C_{j}$ consists of three blocks, each of which
consists of zeros or ones. Writing $0$ for a block of zeros and $1$
for a block of ones, the blocks of $C_{j}$ have either of the two
forms $010$ or $101$, although in the former case the second $0$-block
may be empty, as may be the second $1$-block in the latter case. 

\textbf{$\bullet$ }Two successive columns $C_{j}$ and $C_{j+1}$
cannot both have the $101$-block structure. (This is a consequence
of $n\geq2k$, for since the combined length of the two $1$-blocks
is $k$, the length of the $0$-block in the $101$ case is at least
$k$.) 

When travelling in a common stage $s_{j+1}$, we shall refer to $t_{i}$
and $t_{i'}$ as being \emph{members of the block} if $m_{i,j}$ and
$m_{i',j}$ are in the same block of $C_{j}$. A set of travellers
who are currently moving together will be called a \emph{cohort}.

We now prove inductively on $j$ that, during the period when the
leading cohort is between $P_{j}$ and $P_{j+1}$, the following three
conditions hold:

1. Any pair of members of the same block are in the same cohort. 

2. There are at most $3$ cohorts. 

3. The distance between the members of two neighbouring cohorts is
less than $1$ unit. 

Inductive verification of this trio of claims proves Proposition 5.3. 

For $j=0$ all three claims are clear and indeed there are only $2$
cohorts. Consider $C_{j}$ $(j\geq1)$ and suppose by way of induction
that our claims hold for all lesser values of $j$. 

Take the case where $C_{j}$ has a $010$ block structure. Suppose
first that $C_{j-1}$ also has a $010$ structure. Then, by Lemma
5.2, the members of the joint $01$-block of $C_{j-1}$ arrive together
at $P_{j}$ and from the structure of the transpose cyclic scheme,
form the first $0$-block of $C_{j},$ so forming the lead cohort
in $s_{j+1}$. By induction, the lead of this cohort over the second
$0$-block in $C_{j-1}$ as it becomes the first $0$-block of $C_{j}$
is less than $1$ unit. Since the leading cohort is walking, its lead
over the next cohort remains less than $1$ unit as the lead cohort
traverses $s_{j+1}$. The $1$-block of $C_{j}$ consists of the first
$k$ entries of the members of the second $0$-block of $C_{j-1}$,
whose members arrived in unison at $P_{j}$. This $1$-block cohort
forms the second cohort, which then catches the leading cohort at
$P_{j+1}$. During this time the lead of the second cohort is less
than $1$ unit over the third cohort which is the remainder of the
second $0$-block of $C_{j-1},$ which becomes the second $0$-block
of $C_{j}$ upon arrival at $P_{j}$. These observation taken together
demonstrate that Conditions 1-3 are respected throughout the time
that the lead cohort walks $s_{j+1}$. 

Next suppose that $C_{j-1}$ has a $101$ structure, in which case
the members of the first $1$-block of $C_{j-1}$ arrive first at
$P_{j}$ and form the first $0$-block of $C_{j}$. By induction,
this cohort is less than $1$ unit ahead of the next cohort. Since
the leading block of $C_{j}$ is walking, its lead over the following
cohort cannot increase as it traverses $s_{j+1}$, and so remains
less than $1$ unit. By induction, the separation of the $0$-block
of $C_{j-1}$ and the second $1$-block of $C_{j-1}$ is less than
$1$ unit up until the time the leading cohort of $C_{j-1}$ reaches
$P_{j}$. Their separation decreases after that and the two cohorts
reach $P_{j}$ in unison. After that the $01$-block of $C_{j-1}$
splits into two new cohorts, the first a cohort of size $k$ is comprised
of cyclists, which are the members of the $1$-group of $C_{j}$,
with the remainder of the joint $01$-block of $C_{j-1}$ becoming
the second $0$-block of $C_{j}$ and the third cohort, (thus maintaining
Conditions 1 and 2). The second $0$-block of $C_{j}$ will be less
than $1$ unit behind the second cohort until the leading cohort completes
$s_{j+1}$. Hence Conditions 1, 2, and 3 remain valid throughout the
period where the leading cohort is travelling between $P_{j}$ and
$P_{j+1},$ thus continuing the induction. 

Finally we examine the case where $C_{j}$ has the block form $101$.
By the third bullet point, $C_{j-1}$ has the block form $010$. By
Lemma 5.2, the members of the $01$-block of $C_{j-1}$ arrive together
at $P_{j}$ , and by induction, the members of the second $0$-block
are the trailing cohort, which is less than $1$ unit behind. The
first $1$-block of $C_{j}$ is an initial segment of the joint $01$-block
of $C_{j-1}$ and its members therefore proceed together as the lead
cohort. By construction of the transpose cyclic scheme, the joint
$01$-block of $C_{j-1}$ becomes the joint $10$-block of $C_{j}$,
with the walking members becoming a second cohort in $s_{j+1}$. Their
distance behind the first cohort is always less than $1$ unit. The
second $0$-block of $C_{j-1}$ becomes the second $1$-block of $C_{j}$,
and so its members proceed together, as the third cohort. This cohort
was also the third cohort of $C_{j-1}$ and so was less than $1$
unit behind the members of the $01$-block of $C_{j-1}$ (by Condition
3 and induction) when the joint block reached $P_{j}$. Hence the
separation between the two trailing cohorts is less than $1$ unit
(and decreases to $0$ as these cohorts traverse $s_{j+1}$). Therefore
Conditions 1, 2, and $3$ have been met, and so the induction continues,
thereby completing the proof. $\Box$

\vspace{.2cm}

\textbf{Example 5.4 }Proposition 5.3 does not hold however when $2k>n$.
For example, consider the transpose matrix $M$ for the $(n,n-1)$
problem. Then the zeros consist of the non-leading diagonal running
between entries $(n-1,0)$ and $(0,n-1)$. If we let the ratio of
the cycling speed to walking speed become arbitrarily large, then
$t_{0}$ will reach $P_{n-1}$ before $t_{n-1}$ has reached $P_{1}$,
so that the separation of $t_{0}$ and $t_{n-1}$ approaches an upper
limit of $n-1$ units. 

\section{Calculating features of cyclic schemes }

\textbf{Proposition 6.1 }In respect to the $(n,k)$-cyclic solution,
let $n=r+qk$, $(0\leq r\leq k-1)$, and let $d=$ gcd$(n,k)$. Let
$i_{0},i_{1},\cdots,i_{k-1}$ be the subscripts of the $k$ travellers
that ride stage $s_{1}$. Label the $k$ bicycles as $b_{0},b_{1},\cdots,b_{k-1}$,
where $b_{m}$ is the bicycle ridden by $t_{i_{m}}$ in $s_{1}$.
Then during the execution of the scheme:

(i) the total number of bicycle rides is $n+k-d$.

(ii) Each bicycle $b_{m}$ is mounted on either $\lceil\frac{n}{k}\rceil$
or $\lceil\frac{n}{k}\rceil+1$ occasions, with the first alternative
applying if and only if $r\leq c_{m}$, where $ki_{m}\equiv c_{m}$
(mod $n$), $1\leq c_{m}\leq k$. 

\vspace{.2cm}

\textbf{Proof} (i) There are $k$ travellers $t_{i}$ that cycle $s_{1}$,
and $t_{i}$ completes their quota if and only if $ki\equiv0$ (mod
$n)$. There are $d$ solutions $i$ to this congruence. Therefore
$n-k+d$ travellers have a single ride while $k-d$ travellers have
two separate rides, one beginning and the other ending their journey.
The total number of cycle rides is therefore $n-k+d+2(k-d)=n+k-d.$ 

(ii) Bicycle $b_{m}$ $(0\leq m\leq k-1)$ is mounted at $P_{0}$
by $t_{i_{m}}$ who dismounts at $P_{c_{m}}$, where $ki_{m}\equiv c_{m}$
(mod $n$) ($1\leq c_{m}\leq k$). If $c_{m}<r$ then 
\begin{equation}
c_{m}+qk<r+qk=n.
\end{equation}
Hence $b_{m}$ is ridden by $1+q+1=q+2$ travellers. Since $1\leq r$
we have $\lceil\frac{n}{k}\rceil+1=(q+1)+1=q+2$, as required. Otherwise
$r\leq c_{m}$ whence 
\begin{equation}
c_{m}+(q-1)k\leq qk\leq r+qk=n.
\end{equation}
If $1\leq r$ it follows from $(27)$ that $b_{m}$ is ridden by $1+(q-1)+1=q+1$
travellers. If $r=0$ then $c_{m}=k$ and this figure is $1+(q-1)=q$,
but in either event this number equals $\lceil\frac{n}{k}\rceil$,
thus completing the proof. $\Box$

\vspace{.2cm}

\textbf{Proposition 6.2 }For the $M=M_{n,k}^{T}$ transpose cyclic
matrix and scheme:

(i) If $n\leq2k$ then traveller $t_{i}$ has $k$ cycle rides; if
$2k\geq n$, then $t_{i}$ has $n-k+1$ rides if $n-k\leq i\leq k-1$,
and $n-k$ rides otherwise. 

(ii) The total number of cycle rides is $nk$ if $2k\leq n$ and is
$n(n-k-1)+2k$ if $2k\geq n$. 

(iii) The number of excess handovers $h(M)=\text{min\ensuremath{(}}k(k-1),(n-k)(n-k-1))$. 

(iv) The number of cycle rides after the elimination of excess handovers
is, in all cases, $k(n-k+1)$. 

\vspace{.2cm}

\textbf{Proof }(i) For the case where $2k\leq n$, each traveller
rides just one stage at a time, and so $t_{i}$ has $k$ rides. We
analyse this case further. Applying Proposition 4.3(iv) to $M_{n,k}^{T}$,
we have that $m_{i,j}=m_{n-1-i,n-1-j}$, and so 
\[
m_{0,0}=m_{1,0}=\cdots=m_{k-2,0}=m_{k-1,0}=1,
\]
\begin{equation}
m_{n-k,n-1}=m_{n-k+1,n-1}=\cdots=m_{n-2,n-1}=m_{n-1,n-1}=1.
\end{equation}
\begin{equation}
m_{k,0}=m_{k+1,0}=\cdots=m_{n-1,0}=0=m_{0,n-1}=m_{1,n-1}=\cdots=m_{n-k-1,n-1.}
\end{equation}
It follows from (28) and (29) that for $0\leq i\leq k-1$, row $R_{i}$
has $m_{i,0}=1$, and that $R_{i}$ has $k$ (non-consecutive) entries
equal to $1$, each followed by a maximal sequence of positive length
that consists of entries that equal $0$. For $n-k\leq i\leq n-1$,
the same is true for $R_{i}$ but the statement applies for $R_{i}$
considered in reverse order, beginning with $m_{i,n-1}=1$. On the
other hand, for $k\leq i\leq n-k-1$, $R_{i}$ begins and ends with
a sequence of zeros, and once again there are $k$ entries equal to
$1$, but with no two consecutive entries equal to $1$. 

If we now pass from $M_{n,k}^{T}$ to $\overline{M_{n,k}^{T}}$, the
rows indexed by $0\leq i\leq k-1$ and $n-k\leq i\leq n-1$ each indicate
$k$ bicycle rides, while the remaining central rows each show $k+1$
cycle rides. By symmetry, the same conclusion applies to the matrix
with rows reversed. Now by Remark 4.5 and Proposition 4.3(iii) we
infer that 
\begin{equation}
M_{n,k}^{T}=\overline{(M_{n,n-k})_{r}}^{T}=\overline{((M_{n,n-k})_{r})^{T}}=\overline{((M_{n,n-k}^{T})_{r}}=(\overline{(M_{n,n-k}^{T})})_{r}.
\end{equation}
Therefore if $n<2k$ we infer that the rows $R_{i}$ of $M_{n,k}^{T}$
such that $0\leq i\leq n-k-1$ or $k\leq i\leq n-1$ indicate $n-k$
bicycle rides, while those indexed by $n-k\leq i\leq k-1$ show $n-k+1$
cycle rides, as required.

(ii) If $2k\leq n$ then each cycle is mounted on $n$ separate occasions,
and so the total number of cycle rides is $nk$. Otherwise it follows
by (i) that the total number of cycle rides is given by:
\[
n(n-k)+(k-1-(n-k-1))=(n-k)(n-1)+k=n(n-k-1)+2k.
\]
(iii) The cyclic sequence of $1$'s in $C_{j}$ has length $k$ and
may be written as $T_{j}=(i_{j},i_{j}+1,\cdots,i_{j}+k-1)$ with addition
mod $n$ ($i_{j}=jk$ (mod $n)$). Suppose that $n\geq2k$, in which
case $T_{j}=X_{1,0}^{j}$. By Lemma 5.2 it follows that $s_{j}$ has
no excess handovers unless for some $t$ such that $0\leq t\leq k-1$
we have $i_{j}+2k-1=n+t.$ In that case $j\leq n-2$ and $T_{j+1}$
consists of two linear sequences, which are $I_{1}=(0,1,\cdots,t)$
and $I_{2}=(n-k+t+1,n-k+t+1,\cdots,n-1)$, although the latter is
empty if $t=k-1$. (Note that $|I_{1}|=t+1,$ $|I_{2}|=n-1-(n-k+t)=k-t-1$,
so that $|I_{1}|+|I_{2}|=k$.) Observe from Lemma 5.2 that for all
$i_{r}\in I_{1}$, $i_{r}\in X_{0,1}^{j}$ and $S_{i_{r},j}=S_{i_{j},j}$.
Similarly for all $i_{r}\in I_{2}$, $i_{r}\in X_{1,0}^{j+1}$, and
$i_{s}\in X_{0,1}^{j+1}$ for all $i_{s}\in(t+1,t+2,\cdots,k-1)$,
an interval also of length $k-t-1$, with $S_{i_{r},j+1}=S_{i_{s},j+1}$.
There are no other unnecessary handovers in either $C_{j}$ or in
$C_{j+1}$. It follows that the number of excess handovers in the
pair of stages $s_{j}$ and $s_{j+1}$ is then $|I_{1}|+|I_{2}|=k$.
Conversely if $T_{j+1}$ consists of two linear sequences as above,
then $T_{j}$ is a single linear sequence and the pair $s_{j}$ and
$s_{j+1}$ collectively have $k$ excess handovers. 

It follows that the set of excess handovers of $S(M)$ is partitioned
into sets of order $k$, with one such set for every $0\leq j\leq n-1$
such that $0\leq jk$ (mod $n)\leq k-1$, with one exception. In the
case where $j=n-1$ so that $j+1\equiv0$ (mod $n$) there are not
handovers from $s_{n-1}$ to $s_{0}$. Let $d=$ gcd$(n,k)$. The
number of multiples of $d$ in the interval $[0,k-1]$ is $\frac{k}{d}$.
Then there are $d$ values of $j$ $(1\leq j\leq n-1)$ such that
$jk\equiv td$ (mod $n$) $(0\leq t\leq\frac{n}{d})$. Hence the number
of sets in question is $d\frac{k}{d}=k$. Therefore $h(M)=k^{2}-k=k(k-1)$,
as we subtract $k$ in recognition of no handovers occurring from
$s_{n-1}$ to $s_{0}$. 

On the other hand, if $n\leq2k$ consider $\overline{M_{n,k}^{T}}=(M_{n,n-k}^{T})_{r}$
by (30). Since the latter matrix is also a reverse transpose matrix
and $n\geq2(n-k)$, it now follows from Proposition 3.20 that 
\begin{equation}
h(M_{n,k}^{T})=h(\overline{M_{n,k}^{T}})=h(M_{n,n-k}^{T})_{r}=h(M_{n,n-k}^{T})=(n-k)(n-k-1).
\end{equation}
Therefore $h(M)=k(k-1)$ if $n\geq2k$ and $h(M)=(n-k)(n-k-1)$ otherwise.
Combining these two cases we obtain the statement of (iii). 

For part (iv), if $n\geq2k$ we have the number of cycle rides after
elimination of excess handovers is $kn-k(k-1)=k(n-k+1)$, as required.
In the alternative case, by (i), the corresponding number has the
same form:
\[
n(n-k-1)+2k-(n-k)(n-k-1)
\]
\[
=(n-k-1)(n-n+k)+2k=k(n-k+1).\,\Box
\]

\textbf{Example 6.3 }$M_{11,7}^{T}$ 

\vspace{.2cm}

~~~~~~~~~~$M'=$~~%
\begin{tabular}{|c|c|c|c|c|c|c|c|c|c|c|c|}
\hline 
 & $P_{0}$ & $P_{1}$ & $P_{2}$ & $P_{3}$ & $P_{4}$ & $P_{5}$ & $P_{6}$ & $P_{7}$ & $P_{8}$ & $P_{9}$ & $P_{10}$\tabularnewline
\hline 
\hline 
$t_{0}$ & $1$ & 1 & 0 & 0 & 1 & 1 & 1 & 1 & 0 & 0 & 1\tabularnewline
\hline 
$t_{1}$ & $1$ & 1 & 0 & 0 & 1 & 1 & 1 & 0 & 1 & 1 & 0\tabularnewline
\hline 
$t_{2}$ & $1$ & 1 & 0 & 0 & 1 & 1 & 1 & 0 & 0 & 1 & 1\tabularnewline
\hline 
$t_{3}$ & $1$ & 0 & 1 & 1 & 1 & 0 & 0 & 1 & 1 & 1 & 0\tabularnewline
\hline 
$t_{4}$ & $1$ & 0 & 1 & 1 & 1 & 0 & 0 & 1 & 1 & 1 & 0\tabularnewline
\hline 
$t_{5}$ & $1$ & 0 & 1 & 1 & 0 & 0 & 1 & 1 & 1 & 0 & 1\tabularnewline
\hline 
$t_{6}$ & $1$ & 0 & 0 & 1 & 1 & 1 & 1 & 0 & 0 & 1 & 1\tabularnewline
\hline 
$t_{7}$ & $0$ & 1 & 1 & 1 & 0 & 0 & 1 & 1 & 1 & 0 & 1\tabularnewline
\hline 
$t_{8}$ & $0$ & 1 & 1 & 1 & 0 & 1 & 1 & 0 & 0 & 1 & 1\tabularnewline
\hline 
$t_{9}$ & $0$ & 1 & 1 & 1 & 0 & 1 & 0 & 1 & 1 & 0 & 1\tabularnewline
\hline 
$t_{10}$ & $0$ & 1 & 1 & 0 & 1 & 1 & 0 & 1 & 1 & 1 & 0\tabularnewline
\hline 
\end{tabular}

\vspace{.2cm}

We expunge all excess handovers from $S(M_{11,7}^{T})$ to yield $M'$.
Since $2k\geq n$, from Proposition 6.2(iv) we find that the total
ride number is $11(11-7-1)+2(7)=47$. The number of rides by $t_{0}$
through to $t_{10}$ is $(4+4+4+4)+(5+5+5)+(4+4+4+4)=47$, in accord
with part (i), as $n-k=4$ and $n-k\leq i\leq k-1$ becomes $4\leq i\leq6,$
so it is $t_{4},t_{5},$ and $t_{6}$ who have the extra ride. We
have $h(M)=(11-7)(11-7-1)=12$. According to (iii), after elimination
of excess handovers the total number of rides is $7(11-7+1)=35$,
which indeed equals $47-h(M)$. All travellers have three rides in
$S(M')$ except for $t_{5}$ and $t_{9}$ who each have four. 

\vspace{.2cm}

Throughout this paper we have placed the staging posts at intervals
of one unit with the journey regarded as being of length $n$. We
may however consider other partitions of the travellers' journey.
Consider a putative scheme, $S=S_{m}$, based on partitions into $m$
equal stages. Such a scheme $S_{m}$ would then be represented by
an $n\times m$ binary matrix $M=M(S_{m})$.

\vspace{.3cm}

\textbf{Theorem 6.4 }For the $(n,k)$-problem, let $k=$ gcd$(n,k)$
and put $n'=\frac{n}{d}$ and $k'=\frac{k}{d}$. Then an optimal scheme
$S_{m}$ defined by an $n\times m$ matrix exists for the $(n,k)$-problem
if and only if $m=rn'$ for some $r\geq1$, in which case each traveller
cycles for $l=rk'$ of the $m$ stages of $S_{m}$. 

\vspace{.3cm}

\textbf{Proof} As in the $m=n$ case, for $S_{m}$ to be an optimal
solution, we must have each column $C_{j}$ of $M$ containing exactly
$k$ instances of $1$, and each row containing a common number, $t$
say, of $1$'s. Counting the $1$'s by rows, and then by columns we
equate to see that $tn=km$, whence $m=\frac{tn}{k}=\frac{tn'}{k'}.$
Since gcd$(n',k')=1$, it follows that $k'|t$ so that $t=rk'$ say,
and $m$ necessarily has the form $m=rn'$, for some $r\geq1$. Moreover,
in any optimal scheme $S_{m}$, each traveller cycles the same number,
$l$ say, of stages of $S_{m}$. By optimality we then have $\frac{l}{m}=\frac{k}{n}=\frac{k'}{n'}$
so that $l=\frac{mk'}{n'}=\frac{rn'k'}{n'}=rk'.$ In conclusion:
\begin{equation}
m=rn',\,l=rk'\:(r\geq1).
\end{equation}
Conversely we now show that given that $m$ satisfies (32), we may
build a scheme $S_{m}$ from copies of schemes for the $(n',k')$-problem
to yield an optimal solution for the $(n,k)$-problem based on an
$n\times m$ binary matrix $M$ which is $(l,k)$-uniform, meaning
that each row and each column contains exactly $l$ and $k$ non-zero
entries respectively. To do this we take a $d\times r$ array and
at each position in the array we place an optimal $n'\times n'$ matrix
for the $(n',k')$-problem. (There is no need for these matrices to
be identical solutions.) This yields a $(dn'\times rn')=(n\times m)$
binary matrix $M$ with $l$ entries of $1$ in each row, and $k$
entries of $1$ in each column. 

The first set of $n'$ columns represents a scheme for the initial
$\frac{1}{r}$ part of the journey. Executing this partial scheme
will see $d$ (disjoint) sets of travellers, with each set executing
an optimal $(n',k')$ scheme. Since these schemes are carried out
in parallel, all $n$ travellers will complete the first $\frac{1}{r}$
of the full journey simultaneously, as all these schemes are optimal.
Each of these $d$ sets of travellers will then repeat a similar process
for the second and subsequent partial schemes, with all travellers
completing each of the fractional journeys of lengths $\frac{1}{r},\frac{2}{r},\cdots,\frac{j}{r},\cdots\frac{r}{r}=1$
at the same time. All $n$ travellers complete the journey simultaneously,
having cycled $l$ stages, as required to finish in the least time.
$\Box$

\end{document}